\input amstex
\documentstyle{amsppt}
%
%
\nopagenumbers
\hyphenation{dy-na-mi-cal}
\font\titfont=cmss12
\def\negskp{\hskip -2pt}
\def\compos{\,\raise 1pt\hbox{$\sssize\circ$} \,}
\pagewidth{360pt}
\pageheight{606pt}
\rightheadtext{A note on dynamical systems \dots}
\topmatter
\title\nofrills
\titfont
A note on Newtonian, Lagrangian, and Hamiltonian
dynamical systems in Riemannian manifolds.
\endtitle
\author
R.~A.~Sharipov
\endauthor
\abstract
Newtonian, Lagrangian, and Hamiltonian dynamical systems are
well formalized mathematically. They give rise to geometric
structures describing motion of a point in smooth manifolds.
Riemannian metric is a different geometric structure
formalizing concepts of length and angle. The interplay of
Riemannian metric and its metric connection with mechanical
structures produces some features which are absent in the case
of general (non-Riemannian) manifolds. The aim of present
paper is to discuss these features and develop special
language for describing Newtonian, Lagrangian, and Hamiltonian
dynamical systems in Riemannian manifolds.
\endabstract
\address Rabochaya street 5, 450003, Ufa, Russia
\endaddress
\email \vtop to 20pt{\hsize=280pt\noindent
R\_\hskip 1pt Sharipov\@ic.bashedu.ru\newline
r-sharipov\@mail.ru\vss}
\endemail
\urladdr
http:/\negskp/www.geocities.com/CapeCanaveral/Lab/5341
\endurladdr
\endtopmatter
\loadbold
\TagsOnRight
\document
\head
1. Force field of Newtonian dynamical system.
\endhead
    The primary and most transparent way of describing real mechanical
systems is based on Newton laws. Newton's second law yields differential
equation for the motion of small particle with mass $m$ under the action
of force $\bold F$:
$$
\hskip -2em
m\cdot\ddot\bold r=\bold F(\bold r,\dot\bold r).
\tag1.1
$$
Here $\bold r=\bold r(t)$ is a vector of three-dimensional geometric
space marking position of moving particle. Formally, one can consider
the equation \thetag{1.1} for $\bold r\in\Bbb R^n$ and can take $m=1$
for the sake of simplicity. Further one can replace $\Bbb R^n$ by
arbitrary smooth manifold $M$ and write the equation \thetag{1.1} in
local coordinates $x^1,\,\ldots,\,x^n$:
$$
\hskip -2em
\ddot x^k=\Phi^k(x^1,\ldots,x^n,\,\dot x^1,\ldots,\dot x^n),\qquad
k=1,\,\ldots,\,n.
\tag1.2
$$
Once the equations \thetag{1.2} are written, we meet the problem if
interpreting these equations. If $x^1,\,\ldots,\,x^n$ are coordinates
of moving point $p=p(t)$ in $M$, then their first derivatives $\dot
x^1,\,\ldots,\,\dot x^n$ are components of velocity vector $\bold v
\in T_p(M)$. But second derivatives $\ddot x^1,\,\ldots,\,\ddot x^n$
are not components of a tangent vector of $T_p(M)$. Therefore we are
to consider the pair $q=(p,\bold v)$ being a point of tangent bundle
$TM$, and then write the equations \thetag{1.2} as a system
of first order ODE's:
$$
\hskip -2em
\aligned
&\dot x^k=v^k,\\
&\dot v^k=\Phi^k(x^1,\ldots,x^n,\,v^1,\ldots,v^n).
\endaligned
\tag1.3
$$
Ordinary differential equations \thetag{1.3} correspond to the
following vector field in $TM$:
$$
\pagebreak
\hskip -2em
\boldsymbol\Phi=v^1\cdot\frac{\partial}{\partial x^1}+\ldots
+v^n\cdot\frac{\partial}{\partial x^n}+\Phi^1\cdot\frac{\partial}
{\partial v^1}+\ldots+\Phi^n\cdot\frac{\partial}{\partial v^n}.
\tag1.4
$$
If $q=(p,\bold v)$ is a point of tangent bundle $TM$ and if $\pi\!:
TM\to M$ is a map of canonical projection, then, applying associated
linear map $\pi_*\!:T_q(TM)\to T_p(M)$ to the above vector \thetag{1.4},
we obtain the equality
$$
\hskip -2em
\pi_*\boldsymbol\Phi=\bold v.
\tag1.5
$$
\definition{Definition 1.1} Vector field $\boldsymbol\Phi$ in
tangent bundle $TM$ satisfying the condition \thetag{1.5} is
called {\it Newtonian vector field}.
\enddefinition
\definition{Definition 1.2} {\it Newtonian dynamical system}
in smooth manifold $M$ is a dynamical system determined by
some Newtonian vector field in $TM$.
\enddefinition
    For the motion of real particle both vectors $\bold v$ and
$\bold F$ are in the same space. We can measure their lengths
and the angle between them. Passing to general case of
$n$-dimensional smooth manifold $M$, we loose this opportunity.
Indeed, vector $\boldsymbol\Phi$ is $2n$-dimensional vector
tangent to $TM$, while $\bold v$ is $n$-dimensional vector
tangent to $M$. This situation changes crucially if we take
Riemannian manifold $M$. In this case we can consider vector
$\bold F$ with components
$$
\hskip -2ex
F^k=\Phi^k+\sum^n_{i=1}\sum^n_{j=1}\Gamma^k_{ij}\,v^i\,v^j,
\qquad k=1,\,\ldots,\,n.
\tag1.6
$$
It is tangent to $M$ at the point $p=\pi(q)$. But its components
\thetag{1.6} are functions of double set of arguments $x^1,\,
\ldots,\,x^n,\,v^1,\,\ldots,\,v^n$. In other words, $\bold F$ is
a vector in $T_p(M)$ depending on the point $q=(p,\bold v)\in TM$.
\definition{Definition 1.3} {\it Extended} vector field $\bold F$
in $M$ is a vector-valued function that maps each point $q\in G
\subseteq TM$ to a vector of tangent space $T_p(M)$, where
$p=\pi(q)$. Subset $G$ of $TM$ is a domain of extended vector
field $\bold F$. If $G=TM$, then $\bold F$ is called {\it global}
extended vector field in $M$. 
\enddefinition
    Vector $\bold F$ with components \thetag{1.6} is called {\it force
vector}. It determines {\it force field} of Newtonian dynamical system
in Riemannian manifold. Force field $\bold F$ of Newtonian dynamical
system is an extended vector field in the sense of definition~1.3.
Velocity vector $\bold v$ can also be treated as extended vector field.
Indeed, if $q=(p,\bold v)$ is a point of tangent bundle $TM$, then one
can map it to the vector $\bold v\in T_p(M)$. Now we can calculate
modulus of velocity vector $\bold v$ and scalar product of vectors
$\bold v$ and $\bold F$. In terms of force field $\bold F$ differential
equations \thetag{1.3} are written as:
$$
\hskip -2em
\aligned
&\dot x^k=v^k,\\
&\nabla_tv^k=F^k(x^1,\ldots,x^n,\,v^1,\ldots,v^n).
\endaligned
\tag1.7
$$
Here $\nabla_tv^k$ are components of vector $\nabla_t\bold v$, where
$\nabla_t$ is a covariant derivative with respect to time variable
$t$ along trajectory:
$$
\hskip -2em
\nabla_tv^k=\dot v^k+\sum^n_{i=1}\sum^n_{j=1}\Gamma^k_{ij}
\,v^i\,v^j.
\tag1.8
$$
\definition{Definition 1.4} {\it Newtonian dynamical system}
in smooth Riemannian manifold $M$ is a dynamical system determined
by some extended vector \pagebreak field $\bold F$ in $M$.
\enddefinition
     Note that Newtonian dynamical systems \thetag{1.7} in Riemannian
manifolds are not purely artificial objects obtained as mathematical
generalizations of the equation \thetag{1.1}. As shown in
Chapter~\uppercase\expandafter{\romannumeral 2} of thesis \cite{1},
they arise in describing constrained mechanical systems with holonomic
constraints. Riemannian metric in configuration space of such systems
is given by quadratic form of {\it kinetic energy}.
\head
2. Extended tensor fields.
\endhead
    Extended tensor fields are defined in a similar way as extended
vector fields in definition~1.3. Let's denote by $T^r_s(p,M)$ the
following tensor product:
$$
T^r_s(p,M)=\overbrace{T_p(M)\otimes\ldots\otimes
T_p(M)}^{\text{$r$ times}}\otimes\underbrace{T^*_p(M)
\otimes\ldots\otimes T^*_p(M)}_{\text{$s$ times}}
$$
Tensor product $T^r_s(p,M)$ is the space of tensors of type $(r,s)$ at
the point $p\in M$.
\definition{Definition 2.1} {\it Extended} tensor field $\bold X$
of type $(r,s)$ in $M$ is a tensor-valued function that maps each
point $q\in G\subseteq TM$ to a tensor of the space $T^r_s(p,M)$,
where $p=\pi(q)$. Subset $G$ of $TM$ is a domain of extended tensor
field $\bold X$. If $G=TM$, then $\bold X$ is called {\it global}
extended tensor field in $M$. 
\enddefinition
    Traditional tensor fields of type $(r,s)$ in $M$ are sections
of tensor bundle $T^r_sM$. Extended tensor fields of type $(r,s)$
are sections of pull-back tensor bundle $\pi_*(T^r_sM)$ induced
by the map of canonical projection $\pi\!:TM\to M$. Below we recall
some facts concerning extended tensor fields. Detailed explanation
of the theory of such fields can be found in Chapters~\uppercase
\expandafter{\romannumeral 2}, \uppercase\expandafter{\romannumeral
3}, and \uppercase\expandafter{\romannumeral 4} of thesis \cite{1}.
\par
    The most important fact of the theory of extended tensor fields
in Riemannian manifolds is the presence of two covariant
differentiations
$$
\xalignat 2
&\nabla\!:\ T^r_s(M)\to T^r_{s+1}(M),
&&\tilde\nabla\!:\ T^r_s(M)\to T^r_{s+1}(M).
\endxalignat
$$
First covariant differentiation $\nabla$ is called {\it spatial
differentiation} or {\it spatial gradient}. In local coordinates
it is represented by formula
$$
\hskip -2em
\aligned
&\nabla_qX^{i_1\ldots\,i_r}_{j_1\ldots\,j_s}=\frac{\partial
X^{i_1\ldots\,i_r}_{j_1\ldots\,j_s}}{\partial x^q}
-\sum^n_{a=1}\sum^n_{b=1}v^a\,\Gamma^b_{qa}\,\frac{\partial
X^{i_1\ldots\,i_r}_{j_1\ldots\,j_s}}{\partial v^b}\,+\\
&+\sum^r_{k=1}\sum^n_{a_k=1}\!\Gamma^{i_k}_{q\,a_k}\,X^{i_1\ldots\,
a_k\ldots\,i_r}_{j_1\ldots\,\ldots\,\ldots\,j_s}
-\sum^s_{k=1}\sum^n_{b_k=1}\!\Gamma^{b_k}_{q\,j_k}
X^{i_1\ldots\,\ldots\,\ldots\,i_r}_{j_1\ldots\,b_k\ldots\,j_s}.
\endaligned
\tag2.1
$$
Second covariant differentiation $\tilde\nabla$ is given by much
more simple formula:
$$
\hskip -2em
\tilde\nabla_qX^{i_1\ldots\,i_r}_{j_1\ldots\,j_s}=\frac{\partial
X^{i_1\ldots\,i_r}_{j_1\ldots\,j_s}}{\partial v^q}.
\tag2.2
$$
It is called {\it velocity differentiation} or {\it velocity gradient}.
Velocity gradient $\tilde\nabla$ is defined in arbitrary smooth manifold.
Unlike $\nabla$, it doesn't require the presence of Riemannian metric
in the manifold.\par
\head
3. Covariant representation of extended tensor fields.
\endhead
    Note, that if we replace tangent bundle $TM$ by cotangent
bundle $T^*\!M$, we obtain another definition of extended tensor
fields in $M$.
\definition{Definition 3.1} {\it Extended} tensor field $\bold X$
of type $(r,s)$ in $M$ is a tensor-valued function that maps each
point $q\in G\subseteq T^*\!M$ to a tensor of the space $T^r_s(p,M)$,
where $p=\pi(q)$. Subset $G$ of $T^*\!M$ is a domain of extended
tensor field $\bold X$. If $G=T^*\!M$, then $\bold X$ is called
{\it global} extended tensor field in $M$. 
\enddefinition
    In the case of arbitrary smooth manifold $M$ definitions~2.1 and
3.1 lead to different theories. But for Riemannian manifold $M$ tangent
bundle $TM$ and cotangent bundle $T^*\!M$ are bound with each other
by duality maps
$$
\xalignat 2
&\hskip -2em
\bold g\!:TM\to T^*\!M,
&&\bold g^{-1}\!:T^*\!M\to TM.
\tag3.1
\endxalignat
$$
In local coordinates duality maps \thetag{3.1} are represented as
index lowering and index raising procedures applied to the components
of velocity vector $\bold v$:
$$
\xalignat 2
&v_a=\sum^n_{c=1}g_{ac}\,v^c,
&&v^c=\sum^n_{a=1}g^{ca}\,v_a.
\endxalignat
$$
Due to duality maps \thetag{3.1} two objects introduced by
definitions~2.1 and 2.2 are the same in essential. We call
them {\it contravariant} and {\it covariant} representations
of extended tensor field $\bold X$. If we take covariant
representation of $\bold X$, then formula \thetag{2.1} for
spatial covariant differentiation $\nabla$ is rewritten as
$$
\hskip -2em
\aligned
&\nabla_qX^{i_1\ldots\,i_r}_{j_1\ldots\,j_s}=\frac{\partial
X^{i_1\ldots\,i_r}_{j_1\ldots\,j_s}}{\partial x^q}
+\sum^n_{a=1}\sum^n_{b=1}v_a\,\Gamma^a_{qb}\,\frac{\partial
X^{i_1\ldots\,i_r}_{j_1\ldots\,j_s}}{\partial v_b}\,+\\
&+\sum^r_{k=1}\sum^n_{a_k=1}\!\Gamma^{i_k}_{q\,a_k}\,X^{i_1\ldots\,
a_k\ldots\,i_r}_{j_1\ldots\,\ldots\,\ldots\,j_s}
-\sum^s_{k=1}\sum^n_{b_k=1}\!\Gamma^{b_k}_{q\,j_k}\,
X^{i_1\ldots\,\ldots\,\ldots\,i_r}_{j_1\ldots\,b_k\ldots\,j_s}.
\endaligned
\tag3.2
$$
Formula \thetag{2.2} for velocity gradient $\tilde\nabla$
now is written as follows:
$$
\hskip -2em
\tilde\nabla_qX^{i_1\ldots\,i_r}_{j_1\ldots\,j_s}=
\sum^n_{k=1}g_{qk}\,
\frac{\partial X^{i_1\ldots\,i_r}_{j_1\ldots\,j_s}}
{\partial v_k}.
\tag3.3
$$
In order to make formulas \thetag{2.2} and \thetag{3.3} more
similar to each other we raise index $q$ in \thetag{3.3}. As
a result we get the following formula for $\tilde\nabla$:
$$
\hskip -2em
\tilde\nabla^qX^{i_1\ldots\,i_r}_{j_1\ldots\,j_s}=
\frac{\partial X^{i_1\ldots\,i_r}_{j_1\ldots\,j_s}}
{\partial v_q}.
\tag3.4
$$
\head
4. Differentiation along curves.
\endhead
    Let $p=p(t)$ be some parametric curve (e\.\,g\. the trajectory of
Newtonian dynamical system \thetag{1.7}). Suppose that at each point
$p(t)$ of this curve some tensor $\bold X=\bold X(t)$ of type $(r,s)$
is given. If $\bold X(t)$ is smooth function of $t$, then one can
differentiate it with respect to parameter $t$ along the curve. This
is done by means of covariant derivative $\nabla_t$. As a result we
get another tensor-valued function $\nabla_t\bold X$ on the curve.
Its components are given by the following well-known formula:
$$
\hskip -2em
\gathered
\nabla_tX^{i_1\ldots\,i_r}_{j_1\ldots\,j_s}=\frac{dX^{i_1\ldots\,
i_r}_{j_1\ldots\,j_s}}{dt}+\sum^n_{q=1}\sum^r_{k=1}\sum^n_{a_k=1}
\!\Gamma^{i_k}_{q\,a_k}\,X^{i_1\ldots\,a_k\ldots\,i_r}_{j_1\ldots
\,\ldots\,\ldots\,j_s}\,\dot x^q\,-\\
-\,\sum^n_{q=1}\sum^s_{k=1}\sum^n_{b_k=1}\!\Gamma^{b_k}_{q\,j_k}
\,X^{i_1\ldots\,\ldots\,\ldots\,i_r}_{j_1\ldots\,b_k\ldots\,j_s}
\,\dot x^q.
\endgathered
\tag4.1
$$
Formula \thetag{1.8} is a special case of formula \thetag{4.1}
with $\bold X=\bold v(t)$, and with time derivatives $\dot x^q$
being replaced by $v^q$.\par
    Now suppose again that some curve $p=p(t)$ in $M$ is given.
Its tangent vector $\bold v$ with components $\dot x^1,\,\ldots,
\,\dot x^n$ is vector-valued function of parameter $t$. Taking
pairs $q=(p,\bold v)$, where $p=p(t)$ and $\bold v=\bold v(t)$,
we construct a parametric curve $q=q(t)$ in $TM$. This curve is
called {\it natural lift} of initial curve $p=p(t)$.\par
    Suppose that $\bold X$ is some extended tensor field of type
$(r,s)$ in $M$. According to the definition~2.1, it is tensor-valued
function with argument $q\in TM$. Substituting $q=q(t)$ into the
argument of extended tensor field $\bold X(q)$, we get tensor-valued
function $\bold X(t)$. If $q=q(t)$ is natural lift of curve
$p=p(t)$, then tensor-function $\bold X(t)$ is called {\it natural
restriction} of extended tensor field $\bold X$ to the curve
$p=p(t)$. Let's apply $\nabla_t$ to $\bold X(t)$. As a result for
components of tensor field $\nabla_t\bold X$ we get
$$
\hskip -2em
\nabla_tX^{i_1\ldots\,i_r}_{j_1\ldots\,j_s}=
\sum^n_{k=1}\nabla_kX^{i_1\ldots\,i_r}_{j_1\ldots\,j_s}\cdot v^k+
\sum^n_{k=1}\tilde\nabla_kX^{i_1\ldots\,i_r}_{j_1\ldots\,j_s}\cdot
\nabla_tv^k.
\tag4.2
$$
One can easily write formula \thetag{4.2} in coordinate free form.
Here it is:
$$
\hskip -2em
\nabla_t\bold X=C(\nabla\bold X\otimes\bold v)+C(\tilde\nabla\bold X
\otimes\nabla_t\bold v).
\tag4.3
$$
Note that $\nabla$ and $\tilde\nabla$ in right hand sides of formulas
\thetag{4.2} and \thetag{4.3} are spatial and velocity gradients
respectively, while $C$ is the operation of contraction.
\head
5. Lagrangian dynamical systems.
\endhead
   Lagrangian dynamical system is a special case of Newtonian dynamical
system. The equations of dynamics \thetag{1.3} for them are given in
implicit form by equations
$$
\hskip -2em
\frac{d}{dt}\frac{\partial L}{\partial\dot x^k}-\frac{\partial L}
{\partial x^k}=0,\qquad k=1,\,\ldots,\,n.
\tag5.1
$$
Here $L=L(x^1,\ldots,x^n,\dot x^1,\ldots,\dot x^n)$ is Lagrange
function. Differential equations \thetag{5.1} are known as {\it
Euler-Lagrange equations}. Differentiating composite function,
we can rewrite Euler-Lagrange equations as follows
$$
\hskip -2em
\sum^n_{s=1}\frac{\partial^2 L}{\partial v^k\,\partial v^s}
\cdot\dot v^s+\sum^n_{s=1}\frac{\partial^2 L}{\partial v^k\,
\partial x^s}\cdot v^s-\frac{\partial L}{\partial x^k}=0.
\tag5.2
$$
It's clear that $L$ is a scalar function in $TM$. In other
words, it is extended scalar field. Therefore we can rewrite
\thetag{5.2} in terms of covariant differentiations determined
by formulas \thetag{2.1} and \thetag{2.2}. Let's use \thetag{1.8}
for to express $\dot v^s$ in \thetag{5.2} through covariant
derivative $\nabla_tv^s$. Then let's use formula \thetag{2.1}
applied to scalar field $L$ for to express partial derivative
$\partial L/\partial x^k$ through spatial gradient $\nabla_kL$.
This yields
$$
\hskip -2em
\gathered
\sum^n_{s=1}\frac{\partial^2 L}{\partial v^k\,\partial v^s}
\cdot\nabla_tv^s-\sum^n_{i=1}\sum^n_{j=1}\sum^n_{s=1}
\Gamma^s_{ij}\,\frac{\partial^2 L}{\partial v^k\,\partial v^s}
\,v^i\,v^j\,+\\
+\,\sum^n_{s=1}\frac{\partial^2 L}{\partial v^k\,
\partial x^s}\,v^s-\nabla_kL+\sum^n_{s=1}\sum^n_{j=1}
v^j\,\Gamma^s_{kj}\frac{\partial L}{\partial v^j}=0.
\endgathered
\tag5.3
$$
Gathering second, third, and fifth terms in \thetag{5.3} and
using formulas \thetag{2.1} and \thetag{2.2}, we can transform
\thetag{5.3} to the following form:
$$
\hskip -2em
\sum^n_{s=1}\tilde\nabla_s\!\tilde\nabla_kL
\cdot\nabla_tv^s+\sum^n_{s=1}\nabla_s\!\tilde\nabla_kL\cdot v^s
-\nabla_kL=0.
\tag5.4
$$
Now, if we recall formula \thetag{4.2}, we can further simplify
our equations \thetag{5.4}:
$$
\hskip -2em
\nabla_t\bigl(\tilde\nabla_kL\bigr)-\nabla_kL=0.
\tag5.5
$$
This form of Euler-Lagrange equations is quite similar to initial
one. But now these equations are written in terms of covariant
derivatives \thetag{2.1} and \thetag{2.2}.\par
\head
6. Legendre transformation.
\endhead
   In what case Lagrangian dynamical system determined by
Euler-Lagrange equations \thetag{5.5} can be written in
Newtonian form \thetag{1.7}\,? The answer to this question
depend on the value of determinant of the matrix $A$ with
the following components:
$$
\hskip -2em
A_{ij}=\tilde\nabla_i\!\tilde\nabla_jL=\frac{\partial L}
{\partial v^i\,\partial v^j}.
\tag6.1
$$
If $\det A\neq 0$ then, using \thetag{5.4}, we can express
$\nabla_tv^s$ in explicit form, thus obtaining expression
for the force field of corresponding Newtonian dynamical system.
Lagrangian dynamical system for which the condition
$$
\hskip -2em
\det A\neq 0
\tag6.2
$$
is fulfilled is called {\it regular}. Now suppose that $L$ is
Lagrange function for regular Lagrangian dynamical system
\thetag{5.5}. Then let's denote $\bold p=\tilde\nabla L$. It is
clear that $\bold p$ is an extended covector field with components
$$
\hskip -2em
p_k=\tilde\nabla_kL=\frac{\partial L}{\partial v^k},\qquad
1,\,\ldots,\,n.
\tag6.3
$$
Extended covector field $\bold p$ with components \thetag{6.3}
is used to define nonlinear map $\lambda\!: TM\to T^*\!M$. Indeed,
if $q=(p,v)$ is a point of $TM$, then, taking $\bold p=\bold p(q)$,
we can construct another pair $\tilde q=(p,\bold p)$ being a point
of $T^*\!M$:
$$
\hskip -2em
\lambda(q)=\tilde q=(\pi(q),\bold p(q)).
\tag6.4
$$
Nonlinear map $\lambda\!: TM\to T^*\!M$ defined by formula
\thetag{6.4} is called {\it Legendre transformation}. The
above condition \thetag{6.2} provides local invertibility
of Legendre transformation. Traditionally matrix \thetag{6.1}
is assumed to be a positive matrix (see \cite{2}):
$$
\hskip -2em
A>0.
\tag6.5
$$
This means that $A$ is a matrix of positive quadratic form.
Under the condition \thetag{6.5} Legendre transformation
\thetag{6.4} is globally invertible, i\.e\. it is nonlinear
bijective map $T_p(M)\to T^*_p(M)$ at each point $p\in M$.
The whole set of maps binding tangent and cotangent bundles
(including linear maps \thetag{3.1}) is shown on diagram
below:
%
%
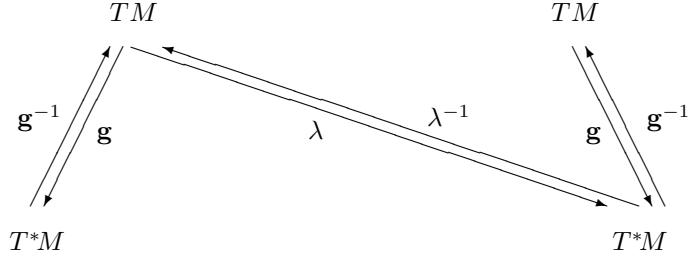
\catcode`@=11
\newcount\@xarg
\newcount\@yarg
\newcount\@yyarg
\newcount\@tempcnta
\newcount\@tempcntb
\newdimen\@linelen
\newdimen\unitlength\unitlength =1pt
\newdimen\@wholewidth\newdimen\@halfwidth
\newdimen\@clnwd
\newdimen\@clnht
\newdimen\@tempdima
\newdimen\@tempdimb
\newdimen\@picht
\newif\if@negarg
\newbox\@linechar
\newbox\@tempboxa
\newbox\@picbox
\font\tenln=line10
\font\tenlnw=linew10
\def\mbox#1{\leavevmode\hbox{#1}}
\def\@ifnextchar#1#2#3{\let\@tempe #1\def\@tempa{#2}\def\@tempb{#3}\futurelet
    \@tempc\@ifnch}
\def\@ifnch{\ifx \@tempc \@sptoken \let\@tempd\@xifnch
      \else \ifx \@tempc \@tempe\let\@tempd\@tempa\else\let\@tempd\@tempb\fi
      \fi \@tempd}
\def\picture(#1,#2){\@ifnextchar({\@picture(#1,#2)}{\@picture(#1,#2)(0,0)}}
\def\@picture(#1,#2)(#3,#4){\@picht #2\unitlength
\setbox\@picbox\hbox to#1\unitlength\bgroup
\hskip -#3\unitlength \lower #4\unitlength \hbox\bgroup\ignorespaces}
\def\endpicture{\egroup\hss\egroup\ht\@picbox\@picht
\dp\@picbox\z@\mbox{\box\@picbox}}
\long\def\put(#1,#2)#3{\@killglue\raise#2\unitlength\hbox to\z@{\kern
#1\unitlength #3\hss}\ignorespaces}
\def\@killglue{\unskip\@whiledim \lastskip >\z@\do{\unskip}}
\def\@whilenoop#1{}
\def\@whiledim#1\do #2{\ifdim #1\relax#2\@iwhiledim{#1\relax#2}\fi}
\def\@iwhiledim#1{\ifdim #1\let\@nextwhile\@iwhiledim
        \else\let\@nextwhile\@whilenoop\fi\@nextwhile{#1}}
\def\thinlines{\let\@linefnt\tenln \let\@circlefnt\tencirc
  \@wholewidth\fontdimen8\tenln \@halfwidth=0.5\@wholewidth}
\def\thicklines{\let\@linefnt\tenlnw \let\@circlefnt\tencircw
  \@wholewidth\fontdimen8\tenlnw \@halfwidth=0.5\@wholewidth}
\def\@getrarrow(#1,#2){\@tempcntb #2\relax
\ifnum\@tempcntb <\z@ \@tempcntb -\@tempcntb\relax\fi
\ifcase \@tempcntb\relax \@tempcnta'55 \or
\ifnum #1<\thr@@ \@tempcnta #1\relax\multiply\@tempcnta
24\advance\@tempcnta -6 \else \ifnum #1=\thr@@ \@tempcnta 49
\else\@tempcnta 58 \fi\fi\or
\ifnum #1<\thr@@ \@tempcnta=#1\relax\multiply\@tempcnta
24\advance\@tempcnta -\thr@@ \else \@tempcnta 51 \fi\or
\@tempcnta #1\relax\multiply\@tempcnta
\sixt@@n \advance\@tempcnta -\tw@ \else
\@tempcnta #1\relax\multiply\@tempcnta
\sixt@@n \advance\@tempcnta 7 \fi\ifnum #2<\z@ \advance\@tempcnta 64 \fi
\char\@tempcnta}
\def\@getlinechar(#1,#2){\@tempcnta#1\relax\multiply\@tempcnta 8%
\advance\@tempcnta -9\ifnum #2>\z@ \advance\@tempcnta #2\relax\else
\advance\@tempcnta -#2\relax\advance\@tempcnta 64 \fi
\char\@tempcnta}
\def\vector(#1,#2)#3{\@xarg=#1\relax\@yarg=#2\relax
\@tempcnta=\ifnum\@xarg<\z@ -\@xarg\else\@xarg\fi
\ifnum\@tempcnta<5\relax
\@linelen=#3\unitlength
\ifnum\@xarg=\z@\@vvector
  \else\ifnum\@yarg=\z@ \@hvector\else\@svector\fi
\fi
\else\@badlinearg\fi}
\def\straightline(#1,#2)#3{\@xarg #1\relax \@yarg #2\relax
\@linelen #3\unitlength
\ifnum\@xarg =\z@ \@vline
  \else \ifnum\@yarg =\z@ \@hline \else \@sline\fi
\fi}
\def\@hvector{\@hline\hbox to\z@{\@linefnt
\ifnum \@xarg <\z@ \@getlarrow(1,0)\hss\else
    \hss\@getrarrow(1,0)\fi}}
\def\@vvector{\ifnum \@yarg <\z@ \@downvector \else \@upvector \fi}
\def\@upvector{\@upline\setbox\@tempboxa\hbox{\@linefnt\char'66}\raise
     \@linelen \hbox to\z@{\lower \ht\@tempboxa\box\@tempboxa\hss}}
\def\@downvector{\@downline\lower \@linelen
      \hbox to \z@{\@linefnt\char'77\hss}}
\def\@upline{\hbox to \z@{\hskip -\@halfwidth \vrule width\@wholewidth
   height \@linelen depth\z@\hss}}
\def\@downline{\hbox to \z@{\hskip -\@halfwidth \vrule \@width \@wholewidth
   \@height \z@ \@depth \@linelen \hss}}
\def\@svector{\@sline
\@tempcnta\@yarg \ifnum\@tempcnta <\z@ \@tempcnta -\@tempcnta\fi
\ifnum\@tempcnta <5%
  \hskip -\wd\@linechar
  \@upordown\@clnht \hbox{\@linefnt  \if@negarg
  \@getlarrow(\@xarg,\@yyarg)\else \@getrarrow(\@xarg,\@yyarg)\fi}%
\else\@badlinearg\fi}
\def\@getlarrow(#1,#2){\ifnum #2=\z@ \@tempcnta'33 \else
\@tempcnta #1\relax\multiply\@tempcnta \sixt@@n \advance\@tempcnta
-9 \@tempcntb #2\relax\multiply\@tempcntb \tw@
\ifnum \@tempcntb >\z@ \advance\@tempcnta \@tempcntb
\else\advance\@tempcnta -\@tempcntb\advance\@tempcnta 64
\fi\fi\char\@tempcnta}
\def\@hline{\ifnum \@xarg <\z@ \hskip -\@linelen \fi
\vrule height\@halfwidth depth\@halfwidth width \@linelen
\ifnum \@xarg <\z@ \hskip -\@linelen \fi}
\def\@sline{\ifnum\@xarg<\z@ \@negargtrue \@xarg -\@xarg \@yyarg -\@yarg
  \else \@negargfalse \@yyarg \@yarg \fi
\ifnum \@yyarg >\z@ \@tempcnta\@yyarg \else \@tempcnta -\@yyarg \fi
\ifnum\@tempcnta>6 \@badlinearg\@tempcnta\z@ \fi
\ifnum\@xarg>6 \@badlinearg\@xarg \@ne \fi
\setbox\@linechar\hbox{\@linefnt\@getlinechar(\@xarg,\@yyarg)}%
\ifnum \@yarg >\z@ \let\@upordown\raise \@clnht\z@
   \else\let\@upordown\lower \@clnht \ht\@linechar\fi
\@clnwd \wd\@linechar
\if@negarg \hskip -\wd\@linechar \def\@tempa{\hskip -2\wd\@linechar}\else
     \let\@tempa\relax \fi
\@whiledim \@clnwd <\@linelen \do
  {\@upordown\@clnht\copy\@linechar
   \@tempa
   \advance\@clnht \ht\@linechar
   \advance\@clnwd \wd\@linechar}%
\advance\@clnht -\ht\@linechar
\advance\@clnwd -\wd\@linechar
\@tempdima\@linelen\advance\@tempdima -\@clnwd
\@tempdimb\@tempdima\advance\@tempdimb -\wd\@linechar
\if@negarg \hskip -\@tempdimb \else \hskip \@tempdimb \fi
\multiply\@tempdima \@m
\@tempcnta \@tempdima \@tempdima \wd\@linechar \divide\@tempcnta \@tempdima
\@tempdima \ht\@linechar \multiply\@tempdima \@tempcnta
\divide\@tempdima \@m
\advance\@clnht \@tempdima
\ifdim \@linelen <\wd\@linechar
   \hskip \wd\@linechar
  \else\@upordown\@clnht\copy\@linechar\fi}
\thinlines
\catcode`\@=\active
%
%
\centerline{\picture(270,117)
\put(45,97){$TM$}
\put(212,97){$TM$}
\put(10,57){$\bold g^{-1}$}
\put(50,87){\vector(-1,-2){30}}
\put(15,27){\vector(1,2){30}}
\put(40,52){$\bold g$}
\put(248,57){$\bold g^{-1}$}
\put(220,87){\vector(1,-2){30}}
\put(255,27){\vector(-1,2){30}}
\put(225,52){$\bold g$}
\put(53,87){\vector(3,-1){180}}
\put(245,27){\vector(-3,1){180}}
\put(120,52){$\lambda$}
\put(165,57){$\lambda^{-1}$}
\put(7,10){$T^*\!M$}
\put(235,10){$T^*\!M$}
\endpicture}\par
    In section~3 above we agreed to make no difference between
covariant and contravariant representations of extended tensor fields.
This means that we consider $\bold X$ and $\bold X\compos\bold g$
as two forms of the same object. Differentiations $\nabla$ and
$\tilde\nabla$ for these two representations of $\bold X$ are defined
by formulas \thetag{2.1}, \thetag{2.2}, \thetag{3.2}, and \thetag{3.4}.
Differentiations $\nabla$ and $\tilde\nabla$ satisfy the following
identities:
$$
\xalignat 2
&\hskip -2em
\nabla(\bold X\compos\bold g)=(\nabla\bold X)\compos\bold g,
&&\tilde\nabla(\bold X\compos\bold g)=(\tilde\nabla\bold X)
\compos\bold g.
\tag6.6
\endxalignat
$$
Presence of nonlinear maps $\lambda$ and $\lambda^{-1}$ on diagram
above increases the number of representations of extended tensor
field $\bold X$. If $\bold Y=\bold X\compos\lambda^{-1}$, then we
say that $\bold Y$ is $\bold p$-representation or {\it momentum}
representation for $\bold X$, while $\bold X$ is $\bold
v$-representation or velocity representation for $\bold Y$.
For instance, $\nabla L$ is a $\bold v$-representation for covector
field $\bold p$, which is called the field of {\it momentum}. In
$\bold p$-representation components of covector field $\bold p$
are treated as independent variables $p_1,\,\ldots,\,p_n$.\par
    Legendre transformation $\lambda$ does not commute with
differentiations $\nabla$ and $\tilde\nabla$. Unlike \thetag{6.6},
here we have the following equalities:
$$
\align
&\hskip -2em
\tilde\nabla(\bold Y\compos\lambda)=C(\tilde\nabla\bold Y\compos
\lambda\otimes\tilde\nabla\tilde\nabla L),
\tag6.7\\
&\hskip -2em
\nabla(\bold Y\compos\lambda)=\nabla\bold Y\compos\lambda
+C(\tilde\nabla\bold Y\compos\lambda\otimes\nabla\tilde\nabla L).
\tag6.8
\endalign
$$
In local coordinates the equalities \thetag{6.7} and \thetag{6.8}
are written as follows:
$$
\pagebreak 
\align
&\hskip -2em
\tilde\nabla_rX^{i_1\ldots\,i_r}_{j_1\ldots\,j_s}
=\sum^n_{k=1}\tilde\nabla_r\tilde\nabla_kL
\cdot\tilde\nabla^kY^{i_1\ldots\,i_r}_{j_1\ldots\,j_s},
\tag6.9\\
&\hskip -2em
\nabla_rX^{i_1\ldots\,i_r}_{j_1\ldots\,j_s}
=\nabla_rY^{i_1\ldots\,i_r}_{j_1\ldots\,j_s}+
\sum^n_{k=1}\nabla_r\tilde\nabla_kL
\cdot\tilde\nabla^kY^{i_1\ldots\,i_r}_{j_1\ldots\,j_s}.
\tag6.10
\endalign
$$
Here $\bold X=\bold Y\compos\lambda$. Similar to \thetag{4.2},
formulas \thetag{6.9} and \thetag{6.10} express the rule of
differentiation for composite functions. They are proved by
direct calculations.\par
   Further let's define the following extended scalar field
in $\bold v$-representation:
$$
\hskip -2em
h=\sum^n_{k=1}v^k\cdot\tilde\nabla_kL-L.
\tag6.11
$$
Then let's take composition of $h$ with Legendre map $\lambda^{-1}$:
$$
\hskip -2em
H=h\compos\lambda^{-1}.
\tag6.12
$$
As a result we get another extended scalar field $H$. It is called
Hamilton function. Hamilton function is traditionally used in
covariant $\bold p$-representation. This means that its argument
is a point $q=(p,\bold p)$ of cotangent bundle $T^*\!M$:
$$
H=H(x^1,\ldots,x^n,\,p_1,\ldots,p_n).
$$
Due to \thetag{6.12} we have $h=H\compos\lambda$. First let's calculate
$\nabla_ih$ directly, using formula \thetag{6.11}, then let's apply
formula \thetag{6.9} for the same purpose. This yields
$$
\align
&\tilde\nabla_ih=\sum^n_{k=1}(\tilde\nabla_iv^k)\cdot
\tilde\nabla_kL+\sum^n_{k=1}v^k\cdot\tilde\nabla_i\tilde\nabla_kL-
\tilde\nabla_iL=\sum^n_{k=1}v^k\cdot\tilde\nabla_i\tilde\nabla_kL,\\
&\tilde\nabla_ih=\sum^n_{k=1}\tilde\nabla_i\tilde\nabla_kL\cdot
\tilde\nabla^kH.
\endalign
$$
Comparing these two formulas for $\tilde\nabla_ih$ and taking into
account that matrix $A$ with components \thetag{6.1} is non-degenerate,
we obtain
$$
\hskip -2em
v^k=\tilde\nabla^kH,\qquad k=1,\,\ldots,\,n.
\tag6.13
$$
Formula \thetag{6.13} is analogous to \thetag{6.3}, it yields an
explicit expression for inverse Legendre map $\lambda^{-1}$ in
local coordinates. Moreover, this formula means that $\tilde\nabla
H$ is a $\bold p$-representation for vector field of velocity
$\bold v$.\par
    Matrix $A$ with components \thetag{6.1} is Jacobi matrix
for Legendre map $\lambda$. Using \thetag{6.13}, we can calculae
components of Jacobi matrix $B$ for inverse Legendre map:
$$
\hskip -2em
B^{ij}=\tilde\nabla^i\tilde\nabla^jH.
\tag6.14
$$
Matrix $B$ with components \thetag{6.14} inherits properties of
matrix $A$, i\.\,e\. $\det B\neq 0$, and if $A>0$, then $B>0$.
\par
   Let's denote $\bold p$-representation of $L$ by $l$. Then let's
transform \thetag{6.11} to $\bold p$-re\-presentation. Note that
$\bold p$-representation for $\bold v$ is $\tilde\nabla H$ and
$\bold p$-representation for $\tilde\nabla L$ is $\bold p$.
Therefore, combining \thetag{6.3}, \thetag{6.13}, \thetag{6.11},
and \thetag{6.12}, we derive
$$
\hskip -2em
l=\sum^n_{k=1}p_k\cdot\tilde\nabla^kH-H.
\tag6.15
$$
Since $l$ by definition is $\bold p$-representation of $L$, we can write
the equality
$$
\hskip -2em
L=l\compos\lambda.
\tag6.16
$$
Formulas \thetag{6.15} and \thetag{6.16} are quite similar to
\thetag{6.11} and \thetag{6.12}. This reflects symmetry of direct
and inverse Legendre maps.\par
   Let's calculate $\nabla_ih$. First let's do it directly, using
formula \thetag{6.11}. Then let's apply formula \thetag{6.10}.
As a result we get
$$
\align
&\nabla_ih=\sum^n_{k=1}v^k\cdot\nabla_i\tilde\nabla_kL-\nabla_iL,\\
&\nabla_ih=\nabla_iH+\sum^n_{k=1}\nabla_r\tilde\nabla_kL\cdot
\tilde\nabla^kH.
\endalign
$$
Comparing these two formulas for $\nabla_ih$, let's take into account
\thetag{6.13}. This yields
$$
\nabla_iH=-\nabla_iL.
\tag6.17
$$
In coordinate free form the equality \thetag{6.17} is written as
$$
\hskip -2em
\nabla H=-(\nabla L)\compos\lambda^{-1}.
\tag6.18
$$
Formula \thetag{6.18} also reflects the symmetry of direct and inverse
Legendre maps.\par
\head
7. Hamiltonian dynamical systems. 
\endhead
    Legendre transformation is used in order to write Lagrangian
dynamical system \thetag{5.5} in Hamiltonian form. Note that
$v^k$ is time derivative of $x^k$. Therefore the above equations
\thetag{6.13} are the equations of dynamics by themselves:
$$
\hskip -2em
\dot x^k=\tilde\nabla^kH(x^1,\ldots,x^n,\,p_1,\ldots,p_n).
\tag7.1
$$
But they are not complete. In order to complete these equations
\thetag{7.1} we are to calculate time derivatives for $p_1,\,
\ldots,\,p_n$. Let's do it using formulas \thetag{6.3} and
\thetag{4.2}:
$$
\hskip -2em
\nabla_tp_k=\sum^n_{s=1}\tilde\nabla_s\!\tilde\nabla_kL\cdot
\nabla_tv^s+\sum^n_{s=1}\nabla_s\!\tilde\nabla_kL\cdot v^s.
\tag7.2
$$
Comparing \thetag{7.2} and \thetag{5.4} and using formula
\thetag{6.17}, we obtain
$$
\hskip -2em
\nabla_tp_k=-\nabla_kH(x^1,\ldots,x^n,\,p_1,\ldots,p_n).
\tag7.3
$$
Now we see that the equations \thetag{7.1} and \thetag{7.3} form
complete system of ODE's. They are called Hamilton equations.
We gather them into a system:
$$
\xalignat 2
&\hskip -2em
\dot x^k=\tilde\nabla^kH,
&&\nabla_tp_k=-\nabla_kH.
\tag7.4
\endxalignat
$$
Hamiltonian dynamical system given by the equations \thetag{7.4} is
called {\it regular} if matrix $B$ with components \thetag{6.14} is
non-degenerate. \pagebreak Each regular Hamiltonian dynamical system
\thetag{7.4} is locally equivalent to some regular Lagrangian dynamical
system \thetag{5.5}, and vice versa, each regular Lagrangian dynamical
system \thetag{5.5} is locally equivalent to some Hamiltonian dynamical
system \thetag{7.4}. If Legendre map defined by \thetag{6.13} is
bijective, then this equivalence is global.
\head
8. Fiberwise spherically symmetric Lagrangians.
\endhead
    Extended tensor field $\bold X$ is called {\it fiberwise spherically
symmetric} if its components depend only on modulus of velocity vector:
$$
X^{i_1\ldots\,i_r}_{j_1\ldots\,j_s}=X^{i_1\ldots\,i_r}_{j_1\ldots\,
j_s}(x^1,\,\ldots,\,x^n,v)\text{, \ where \ }v=|\bold v|.
$$
This means that $X^{i_1\ldots\,i_r}_{j_1\ldots\,j_s}$ is spherically
symmetric function within each fiber of tangent bundle $TM$ for each
fixed point $p\in M$.
Such fields were considered in Chapter~\uppercase\expandafter{\romannumeral
7} of thesis \cite{1}. Now suppose that Lagrange function $L$ of some
Lagrangian dynamical system \thetag{5.5} is fiberwise spherically symmetric
scalar field. Let's write this field in Newtonian form and let's calculate
its force field $\bold F$. Components of $\bold F$ should be obtained from
the equations \thetag{5.4} written as follows:
$$
\hskip -2em
\sum^n_{s=1}\tilde\nabla_s\!\tilde\nabla_kL
\cdot F^s+\sum^n_{s=1}\nabla_s\!\tilde\nabla_kL\cdot v^s
-\nabla_kL=0.
\tag8.1
$$
Let's calculate covariant derivatives $\tilde\nabla_s\!\tilde\nabla_kL$
and $\nabla_s\!\tilde\nabla_kL$, assuming $L$ to be fiberwise spherically
symmetric. For first order derivative $\tilde\nabla_kL$ we have
$$
\hskip -2em
\tilde\nabla_kL=L'\cdot\frac{v_k}{|\bold v|},
\tag8.2
$$
Here by $L'$ we denote partial derivative of the function
$L(x^1,\ldots,x^n,v)$ with respect to its last argument $v$, which
is interpreted as modulus of velocity vector $\bold v$. Let's apply
covariant derivatives $\nabla_s$ and $\tilde\nabla_s$ to \thetag{8.2}.
This yields
$$
\align
&\hskip -2em
\nabla_s\!\tilde\nabla_kL=\nabla_sL'\cdot\frac{v_k}{|\bold v|},
\tag8.3
\\
&\hskip -2em
\tilde\nabla_s\!\tilde\nabla_kL=L''\cdot\frac{v_s\,v_k}{|\bold v|^2}
+\frac{L'}{|\bold v|}\cdot\left(g_{sk}-\frac{v_s\,v_k}{|\bold v|^2}
\right).
\tag8.4
\endalign
$$
Note that \thetag{8.4} are components of matrix $A$ (see \thetag{6.2}
above). In order to invert this matrix let's consider two operators
of orthogonal projection $\bold Q$ and $\bold P$:
$$
\xalignat 2
&\hskip -2em
Q^i_k=\frac{v^i\,v_k}{|\bold v|^2},
&&P^i_k=\delta^i_k-\frac{v^i\,v_k}{|\bold v|^2}.
\tag8.5
\endxalignat
$$
First of them is a projector to the direction of velocity vector
$\bold v$, second is a projector to hyperplane perpendicular to
$\bold v$. Projection operators $\bold Q$ and $\bold P$ with
components \thetag{8.5} are complementary to each other, this means
that
$$
\hskip -2em
\bold Q+\bold P=\bold 1\quad\text{and}\quad\bold Q\compos\bold P
=\bold P\compos\bold Q=\bold 0.
\tag8.6
$$
Comparing \thetag{8.4} with \thetag{8.5}, we find that
$$
\hskip -2em
A_{sk}=\tilde\nabla_s\!\tilde\nabla_kL=L''\cdot Q_{sk}+
\frac{L'}{|\bold v|}\cdot P_{sk}.
\tag8.7
$$
Matrix $A$ with components \thetag{8.7} is non-degenerate if and only
if $L'\neq 0$ and $L''\neq 0$ simultaneously. In this case matrix
$B=A^{-1}$ has the following components:
$$
\hskip -2em
B^{rk}=\frac{1}{L''}\cdot Q^{rk}+\frac{|\bold v|}{L'}\cdot P^{rk}.
\tag8.8
$$
Now, combining \thetag{8.1}, \thetag{8.3}, and \thetag{8.8}, we derive
formula for components of $\bold F$:
$$
F_r=-\sum^n_{s=1}\left(\frac{\nabla_sL'}{L''}-\frac{\nabla_sL}
{|\bold v|\cdot L''}\right)\cdot\frac{v^s\,v_r}{|\bold v|}
-|\bold v|\sum^n_{s=1}\frac{\nabla_sL}{L'}\cdot\left(\frac{v^s\,v_r}
{|\bold v|^2}-\delta^s_r\right).
\tag8.9
$$
This formula \thetag{8.9} is quite similar to the following one:
$$
\hskip -2em
F_r=-|\bold v|\sum^n_{s=1}\frac{\nabla_sW}{W'}\cdot\left(\frac{2\,\,
v^s\,v_r}{|\bold v|^2}-\delta^s_r\right).
\tag8.10
$$
Here $W$ is some fiberwise spherically symmetric scalar field with
$W'\neq 0$. Force fields of the form \thetag{8.10} arise in the theory
of Newtonian dynamical systems {\it admitting normal shift} (see
Chapter~\uppercase\expandafter{\romannumeral 7} of thesis \cite{1}).
Let's find in which case formulas \thetag{8.9} and \thetag{8.10} do
coincide. This occurs if the following equations hold:
$$
\xalignat 2
&\hskip -2em
\frac{\nabla_sL'}{L''}-\frac{\nabla_sL}
{|\bold v|\cdot L''}+\frac{\nabla_sL}{L'}
=\frac{2\,\nabla_sW}{W'},
&&\frac{\nabla_sL}{L'}=\frac{\nabla_sW}{W'}.
\tag8.11
\endxalignat
$$
We consider nontrivial case, when $\nabla W\neq 0$ and $\nabla L
\neq 0$. In this case equations \thetag{8.11} mean that 
spatial gradients of scalar fields $L$ and $L'$ are collinear.
This occurs if and only if $L'=f(L,v)$, where $f=f(u,v)$ is
some smooth function of two variables. This means that we deal
with a class of functions $L=L(x^1,\ldots,x^n,v)$, each of which
is a solution of partial differential equation
$$
\hskip -2em
\frac{\partial L}{\partial v}=f(L,v)
\tag8.12
$$
with some particular function $f=f(u,v)\neq 0$. Note that we should not
solve the equation \thetag{8.12} for particular function $f$. We should
describe the whole set of solutions for all equations of the form
\thetag{8.12}. This is done by formula
$$
\hskip -2em
L=\beta(C(x^1,\ldots,x^n),v).
\tag8.13
$$
Here $\beta=\beta(u,v)$ is a smooth function of two variables with
$\beta'_v\neq 0$ and $\beta''_{vv}\neq 0$, while $C=C(x^1,\ldots,x^n)$
is a function of spatial variables only, i\.\,e\. this is traditional
(not extended) scalar field in $M$. Let's substitute \thetag{8.13} into
the equation
$$
\frac{\nabla_sL'}{L''}-\frac{\nabla_sL}
{|\bold v|\cdot L''}=\frac{\nabla_sL}{L'}.
$$
derived from \thetag{8.11}. This leads to the following
differential equation for $\beta(u,v)$:
$$
\frac{\beta''_{uv}}{\beta''_{vv}}-\frac{\beta'_u}{v\cdot \beta''_{vv}}
=\frac{\beta'_u}{\beta'_v}.
$$
This equation can be transformed so that it can be further integrated:
$$
\hskip -2em
\frac{\beta''_{uv}}{\beta'_u}=\frac{1}{v}+\frac{\beta''_{vv}}{\beta'_v}.
\tag8.14
$$
Integrating \thetag{8.14} with respect to $v$, we obtain
$$
\hskip -2em
\log(\beta'_u)=\log(v)+\log(\beta'_v)+\log(c)\text{, \ where \ }c=c(u).
\tag8.15
$$
Note that, ultimately, in formula \thetag{8.13} for $L$ we substitute
$u=C(x^1,\ldots,x^n)$, where $C(x^1,\ldots,x^n)$ is arbitrary smooth
function. This means that varying function $c(u)$ in \thetag{8.15},
we do not change class of Lagrange functions $L$. Let's choose
$c(u)=1/u$ for the sake of further convenience. Then we get
$$
u\cdot\beta'_u=v\cdot\beta'_v.
$$
This equation is explicitly integrable. Its general solution is
determined by one arbitrary smooth function of one variable
$\phi=\phi(z)$:
$$
\hskip -2em
\beta(u,v)=\phi(u\cdot v).
\tag8.16
$$
Substituting \thetag{8.16} into \thetag{8.13} and further into
\thetag{8.9}, we obtain
$$
\hskip -2em
F_r=-\sum^n_{s=1}\frac{\nabla_sC}{C}\cdot\left(2\,
v^s\,v_r-|\bold v|^2\cdot\delta^s_r\right).
\tag8.17
$$
If we substitute $C=e^{-f}$, where $f=f(x^1,\ldots,x^n)$, we can
rewrite \thetag{8.17} as
$$
\hskip -2em
F_r=-|\bold v|^2\cdot\nabla_rf+\sum^n_{s=1}2\cdot(\nabla_sf\,v^s)
\cdot v_r
\tag8.18
$$
Newtonian dynamical system \thetag{1.7} with force field $\bold F$
given by formula \thetag{8.18} coincides with geodesic flow of
metric $\tilde g=e^{-2f}\cdot\bold g$, which is conformally equivalent
to basic metric $\bold g$ of Riemannian manifold $M$. Thus we have
proved a theorem.\par
\proclaim{Theorem 8.1} Newtonian dynamical system with force field of
the form \thetag{8.10} possess Lagrangian structure with fiberwise
spherically symmetric Lagrange function $L$ if and only if its force
field is given by formula \thetag{8.18}, which is special case of formula
\thetag{8.10} with $W=v\cdot e^{-f}$, where $f=f(x^1,\ldots,x^n)$.
\endproclaim
\head
9. Inverse problem of Lagrangian dynamics.
\endhead
    Theorem~8.1 is not an ultimate result concerning Lagrangian structures
of Newtonian dynamical systems admitting normal shift. First reason
is that formula \thetag{8.10} does not cover general case (see
Chapter~\uppercase\expandafter{\romannumeral 7} of thesis \cite{1}).
General formula for the force field of Newtonian dynamical system
admitting normal shift of hypersurfaces in Riemannian manifold $M$
with $\dim M\geqslant 3$ looks like
$$
\hskip -2em
F_r=\frac{h(W)}{W'}\cdot\frac{v_r}{|\bold v|}-|\bold v|\sum^n_{s=1}
\frac{\nabla_sW}{W'}\cdot\left(\frac{2\,\,v^s\,v_r}{|\bold v|^2}
-\delta^s_r\right),
\tag9.1
$$
where $h=h(w)$ is an arbitrary smooth function of one variable. And
second reason is that in theorem~8.1 we restrict ourselves to the
case of fiberwise spherically symmetric Lagrange functions.\par
    In order to study general case we should substitute \thetag{9.1}
into \thetag{8.1} and consider \thetag{8.1} as a system of PDE's for
unknown Lagrange function. Problem of determining whether the Newtonian
dynamical system with a given force field $\bold F$ admits Lagrangian
structure (and finding Lagrange function if it admits) is known as {\it
inverse problem of Lagrangian dynamics}. As known to me, this problem
is not solved in general case (see more details and references in
\cite{3--7}). Even for special force fields given by explicit formula
\thetag{9.1} it remains unsolved. Solving this problem for force field
\thetag{9.1} is very important since it would open a way for applying
theory from \cite{1} to the description of wave front dynamics and
to some related problems arising in analysis of partial differential
equations (see \cite{8} and \cite{9}).
\head
10. Resume.
\endhead
    Concept of extended tensor field arisen in \cite{10} and \cite{11},
and used in \cite{1} for describing Newtonian dynamical systems is
applicable to Lagrangian and Hamiltonian dynamical systems in Riemannian
manifolds as well. It gives a method (or a language) for describing these
systems in terms of their configuration space $M$ instead of using
geometric structures in tangent bundle $TM$ (exception is Hamiltonian
dynamical systems in abstract simplectic manifolds, when one cannot
separate configuration space within phase space of dynamical system).
As an example of applying suggested method I consider inverse problem
of Lagrangian dynamics for Newtonian dynamical systems admitting normal
shift, and I give partial solution of this problem in class of fiberwise
spherically symmetric Lagrange functions.


\Refs
\ref\no 1\by Sharipov~R.~A.\book Dynamical systems admitting the normal
shift\publ thesis for the degree of Doctor of Sciences in Russia\publaddr
Ufa\yr 1999\moreref English version of thesis is submitted to Electronic
Archive at LANL\footnotemark, see archive file ArXiv:math.DG/0002202 in the
section of Differential Geometry\footnotemark
\endref
\adjustfootnotemark{-1}
\footnotetext{Electronic Archive at Los Alamos National Laboratory of USA
(LANL). Archive is accessible through Internet 
{\bf http:/\negskp/arXiv.org}, it has mirror site 
{\bf http:/\negskp/ru.arXiv.org} at the Institute for Theoretical and
Experimental Physics (ITEP, Moscow).}
\adjustfootnotemark{+1}
\footnotetext{For the convenience of reader we give direct reference
to archive file. This is the following URL address:
{\bf http:/\negskp/arXiv.org/eprint/math\.DG/0002202}\,.}
\adjustfootnotemark{-2}
\ref\no 2\by Arnold~V.~I.\book Mathematical methods of classical
mechanics\publ Nauka publishers\publaddr Moscow\yr 1979
\endref
\ref\no 3\by Filippov~V.~M., Savchin~V.~M., Shorohov~S.~G.
\paper Variational principles for non-potential operators
\inbook Modern problems in mathematics. Recent achievements
\vol 40\yr 1992\pages 3--176\publ VINITI\publaddr Moscow
\endref
\ref\no 4\by Morandi~G., Ferrario~C., Lo Vecchio G., Marmo~G.,
Rubano C.\paper The inverse problem in the calculus of variations
and the geometry of the tangent bundle\jour Phys\. Reports
\vol 188\yr 1990 \issue 3--4\pages 147--284
\endref
\ref\no 5\by Crampin~M.\paper On the differential geometry of the
Euler--Lagrange equations and the inverse problem of Lagrangian
dynamics\jour Journ\. of Phys\. A.\yr 1981 \vol 14\issue 10
\pages 2567--2575
\endref
\ref\no 6\by Sarlet~W.\paper Contribution to the study of
symmetries, first integrals and inverse problem of the
calculus of variations in theoretical mechanics
\jour Acad\. Analecta\yr 1987\vol 49\issue 1\pages 27--57
\endref
\ref\no 7\by Carinena~W.~F., Lopez~C., Martinez~E.
\paper A geometrical characterization of Lagrangian
second-order differential equations
\jour Inverse Problems\yr 1989\vol 5\issue 5\pages 691--705
\endref
\ref\no 8\by Fedoryuk~M.~V. \paper The equations with fast oscillating
solutions\inbook Summaries of Science and Technology. Modern problems
of Mathematics. Fundamental researches. Vol. 34\yr 1988\publ
VINITI\publaddr Moscow
\endref
\ref\no 9\by Arnold~V.~I.\book Singularities of caustics and wave
fronts\yr 1996\publ Phazis publishers\publaddr Moscow
\endref
\ref\no 10\by Finsler~P.\book \"Uber Kurven and Fl\"achen in algemeinen
Raumen\publ Dissertation\publaddr G\"ottingen\yr 1918
\endref
\ref\no 11\by Cartan~E\.\book Les espaces de Finsler\publ Actualites 79
\publaddr Paris\yr 1934
\endref
\endRefs
\enddocument
\end